\documentclass[graybox]{svmult}

\usepackage{mathptmx}       
\usepackage{helvet}         
\usepackage{courier}        
\usepackage{type1cm}        
\usepackage{makeidx}         
\usepackage{graphicx}        
\usepackage{multicol}        
\usepackage[bottom]{footmisc}

\makeindex             

\usepackage{amsmath, amssymb, tikz, pgfplots, stmaryrd, pgffor}
\usetikzlibrary{calc,plotmarks}
\usepgfplotslibrary{groupplots}

\newcommand{\setR}{ {\mathbb R} }
\newcommand{\opskew}{\operatorname{skew}}
\newcommand{\opdiv}{\operatorname{div}}

\newcommand{\that}{\hat t}
\newcommand{\uhat}{\hat u}

\newcommand{\Ehat}{\hat E}
\newcommand{\Hhat}{\hat H}

\begin{document}

\title*{An Explicit Mapped Tent Pitching Scheme for Maxwell Equations}
\titlerunning{MTP scheme for Maxwell Equations}

\author{Jay Gopalakrishnan
\and
Matthias Hochsteger
\and
Joachim Sch\"oberl
\and
Christoph Wintersteiger
}
\authorrunning{J. Gopalakrishnan, M. Hochsteger, J. Sch\"oberl, C. Wintersteiger}
\institute{
  Jay Gopalakrishnan \at
  Fariborz Maseeh Department of Mathematics \& Statistics,
  Portland State University, PO Box 751, Portland OR 97207-0751, USA,
  \email{gjay@pdx.edu}
  \and
  Matthias Hochsteger, Joachim Sch\"oberl, Christoph Wintersteiger \at
  Institute for Analysis and Scientific Computing,
  Technische Universit\"at Wien, Wiedner Hauptstra\ss e 8-10, 1040
  Wien, Austria,
  \email{matthias.hochsteger@tuwien.ac.at, joachim.schoeberl@tuwien.ac.at, christoph.wintersteiger@tuwien.ac.at}
}

\maketitle

\abstract{\quad We present a new numerical method for solving time
  dependent Maxwell equations, which is also suitable for general
  linear hyperbolic equations. It is based on an unstructured
  partitioning of the spacetime domain into tent-shaped regions that
  respect causality.  Provided that an approximate solution is
  available at the tent bottom, the equation can be locally evolved up
  to the top of the tent. By mapping tents to
  a domain which is a tensor product of a spatial domain with a time
  interval, it is possible to construct a fully 
  explicit scheme that advances the solution through unstructured
  meshes. This work highlights a difficulty that arises when standard
  explicit Runge Kutta schemes are used in this context and proposes
  an alternative structure-aware Taylor time-stepping technique.
  Thus 
  explicit methods are constructed
  that allow variable time steps and local refinements without compromising
  high order accuracy in space and time.
  These Mapped Tent
  Pitching (MTP) schemes lead to highly parallel algorithms, which
  utilize modern computer architectures extremely well.
}

\section{Introduction}
\label{mtp_sat:introduction}
Electromagnetic waves propagate at the speed of light. Thus, the field at
a certain point in space and time depends only on field values within a
dependency cone. A tent pitching method introduces a special
``causal'' 
spacetime mesh that respects this finite speed of propagation. It is
not limited to Maxwell equations, but can be applied to general
hyperbolic equations. A tent pitching method requires a numerical
scheme to discretize the equation on that mesh. Discontinuous Galerkin (DG)
methods are of particular interest since they offer a systematic
avenue to build high order methods. For a given initial condition at
the bottom of a tent, the discrete equations may be solved
within each individual tent, up to the tent top.
The computed solution at the tent top provides initial conditions
for the tents that follow later in time. This method is highly parallel,
since many tents can be solved independently. Methods using such
tent-pitched meshes may be traced back to 
\cite{LowriRoeLeer95, Richt94}. More recent
works~\cite{AbediPetraHaber06, MR05, ASHT00} 
develop  Spacetime DG (SDG) methods within tents 
by formulating local variational
problems, for which linear systems are set up and solved. 
Although these systems are local, the matrix size can grow rapidly with
the polynomial order, especially  in four-dimensional spacetime
tents. In this context it is natural to ask if one can develop
explicit schemes (which usually perform well under low memory bandwidth)
that take advantage of tents. 

A key ingredient to answer this question was presented 
in \cite{GSW17}, where  Mapped Tent Pitching (MTP)
schemes were introduced. The MTP discretization, which proceeds by
mapping tents to a spacetime cylinder,   allows one 
to evolve the solution either implicitly or explicitly
within tents. The  memory  
requirements of the explicit MTP scheme
are limited to what is needed for storing 
the spatial mesh, the solution coefficients at one
time step, and the topology of the tents.

In this work, we show that notwithstanding the above-mentioned
advantages of the explicit MTP scheme, one may lose higher order
convergence if a naive time stepping strategy (involving a standard
explicit Runge-Kutta scheme) is used. We then develop
a new Taylor time-stepping for the local
problems within tents. Despite its simplicity, our numerical
experiments  show that 
it delivers optimal order of
convergence.

\section{Mesh generation by tent pitching} \label{sec:mesh}

\begin{figure}
  \centering
  \begin{tikzpicture}[scale=1.5]
    \tikzstyle{bluefill} = [fill=blue!30,fill opacity=0.5]
    \tikzstyle{redfill} = [fill=red!30,fill opacity=0.5]

    \coordinate (p0) at (0.000,0.000);
    \coordinate (p1) at (1.000,0.500);
    \coordinate (p2) at (2.000,0.000);
    \filldraw[bluefill] (p0) -- (p1) -- (p2);
    \coordinate (p3) at (2.000,0.000);
    \coordinate (p4) at (3.000,0.500);
    \coordinate (p5) at (4.000,0.000);
    \filldraw[bluefill] (p3) -- (p4) -- (p5);
    \coordinate (p6) at (4.000,0.000);
    \coordinate (p7) at (4.800,0.320);
    \coordinate (p8) at (5.440,0.000);
    \filldraw[bluefill] (p6) -- (p7) -- (p8);
    \coordinate (p9) at (5.440,0.000);
    \coordinate (p10) at (5.952,0.205);
    \coordinate (p11) at (6.362,0.000);
    \filldraw[bluefill] (p9) -- (p10) -- (p11);
    \coordinate (p12) at (6.362,0.000);
    \coordinate (p13) at (6.689,0.131);
    \coordinate (p14) at (6.951,0.000);
    \filldraw[bluefill] (p12) -- (p13) -- (p14);
    \coordinate (p15) at (0.000,1.000);
    \draw (p15) -- (p1);
    \coordinate (p16) at (2.000,1.000);
    \draw (p1) -- (p16) -- (p4);
    \coordinate (p17) at (4.000,0.720);
    \draw (p4) -- (p17) -- (p7);
    \coordinate (p18) at (5.440,0.461);
    \draw (p7) -- (p18) -- (p10);
    \coordinate (p19) at (6.362,0.295);
    \draw (p10) -- (p19) -- (p13);
    \coordinate (p20) at (6.951,0.262);
    \draw (p13) -- (p20);
    \coordinate (p21) at (1.000,1.000);
    \draw (p15) -- (p21) -- (p16);
    \coordinate (p22) at (3.000,1.000);
    \draw (p16) -- (p22) -- (p17);
    \coordinate (p23) at (4.800,0.781);
    \draw (p17) -- (p23) -- (p18);
    \coordinate (p24) at (5.952,0.500);
    \filldraw[redfill] (p10) -- (p18) -- (p24) -- (p19);
    \draw[thick] (p15) -- (p21) -- (p16) -- (p22) -- (p17) -- (p23) -- (p18) -- (p10) -- (p19) -- (p13) -- (p20);
    \draw ($0.5*(p22)+0.5*(p17)$) node[circle,radius=2,draw,fill=white,inner sep=0.6] {} -- ++(30:0.7) node[right] {$\tau_i$};
    \draw (p18) -- (p24) -- (p19);
    \coordinate (p25) at (6.689,0.393);
    \draw (p19) -- (p25) -- (p20);
    \coordinate (p26) at (4.000,1.000);
    \draw (p22) -- (p26) -- (p23);
    \coordinate (p27) at (5.440,0.756);
    \draw (p23) -- (p27) -- (p24);
    \coordinate (p28) at (6.362,0.557);
    \draw (p24) -- (p28) -- (p25);
    \coordinate (p29) at (6.951,0.524);
    \draw (p25) -- (p29);
    \coordinate (p30) at (4.800,1.000);
    \draw (p26) -- (p30) -- (p27);
    \coordinate (p31) at (5.952,0.762);
    \draw (p27) -- (p31) -- (p28);
    \coordinate (p32) at (6.689,0.655);
    \draw (p28) -- (p32) -- (p29);
    \coordinate (p33) at (5.440,1.000);
    \draw (p30) -- (p33) -- (p31);
    \coordinate (p34) at (6.362,0.819);
    \draw (p31) -- (p34) -- (p32);
    \coordinate (p35) at (6.951,0.786);
    \draw (p32) -- (p35);
    \coordinate (p36) at (5.952,1.000);
    \draw (p33) -- (p36) -- (p34);
    \coordinate (p37) at (6.689,0.918);
    \draw (p34) -- (p37) -- (p35);
    \coordinate (p38) at (6.362,1.000);
    \draw (p36) -- (p38) -- (p37);
    \coordinate (p39) at (6.951,1.000);
    \draw (p37) -- (p39);
    \coordinate (p40) at (6.689,1.000);
    \draw (p38) -- (p40) -- (p39);
    \draw[densely dotted] (1.000,0) -- (1.000,1.000);
    \draw[densely dotted] (2.000,0) -- (2.000,1.000);
    \draw[densely dotted] (3.000,0) -- (3.000,1.000);
    \draw[densely dotted] (4.000,0) -- (4.000,1.000);
    \draw[densely dotted] (4.800,0) -- (4.800,1.000);
    \draw[densely dotted] (5.440,0) -- (5.440,1.000);
    \draw[densely dotted] (5.952,0) -- (5.952,1.000);
    \draw[densely dotted] (6.362,0) -- (6.362,1.000);
    \draw[densely dotted] (6.689,0) -- (6.689,1.000);
    \draw (6.951,0) -- (6.951,1.000);
    \draw[-latex,thick] (0,0)--(7.2989952,0) node[right] {$x$};
    \draw[-latex,thick] (0,0)--(0,1.3) node[left] {$t$};
  \end{tikzpicture}
  \caption{Tent pitched spacetime mesh for a one-dimensional spatial mesh. \label{mtp_sat:fig1d}}
\end{figure}
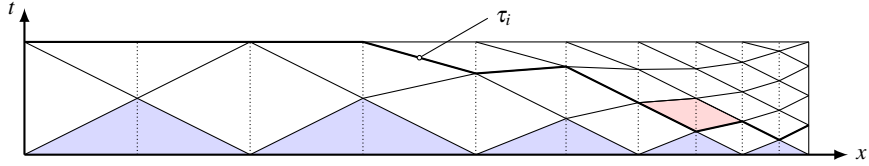

We start with a conforming spatial mesh consisting of 
elements ${\mathcal T} = \{ T \}$ and vertices ${\mathcal V} = \{ V \}$.
We progress in time by defining a sequence of advancing fronts $\tau_i$.
A front $\tau_i$ is given as a standard nodal finite element
function on this mesh. It is defined by storing the current time for every vertex of the
mesh. We move from $\tau_i$ to the next front $\tau_{i+1}$ by moving
one vertex forward in time, while keeping all other vertices
fixed. The spacetime domain between $\tau_{i}$ and $\tau_{i+1}$ we
call a tent. In Fig.~\ref{mtp_sat:fig1d}, the red domain is the
tent between $\tau_i$ and $\tau_{i+1}$.

Its projection to the spatial domain is exactly
the vertex patch $\omega_V$ around $V$ of the original mesh.
The data to be stored for one tent are the bottom and top-times of the
central vertex, plus the times for all neighboring vertices.

Note that although the algorithm is described sequentially, it is
highly parallel. Vertices with graph-distance of at least two can be
moved forward independently. For example, 
in Fig.~\ref{mtp_sat:fig1d}, all blue
tents can be built and processed in parallel.

The distance for advancing a vertex is limited by the speed of light,
a constraint often  referred to in the literature
as the {\em causality condition.}
Under this condition, 
the Maxwell problem inside the tent is solvable using
the initial
conditions at the tent bottom. Thus, the top boundary is an outgoing
boundary and no boundary conditions are needed there.

Note that the spatial mesh is refined towards the right boundary, which leads 
to smaller tent heights at the right boundary. Hence, smaller time steps
in locally refined regions is a very natural feature of tent pitching methods.

\section{The MTP discretization}

Now, we consider the discretization method for one tent domain 
$K = \{ (x,t) : x \in \omega_{V},  \varphi_b(x) \leq t \leq \varphi_t(x) ) \},
$
where $\omega_V$ is the union of elements containing the vertex $V$, and
$\varphi_b$ and $\varphi_t$ are the bottom and top fronts, respectively,
restricted
to $\omega_V$.
Our aim is to numerically solve the Maxwell system on $K$, namely 
\begin{equation} \label{mtp_sat:maxwell}
\partial_t \varepsilon  E = \nabla \times H\;, \qquad \partial_t \mu   H = -\nabla
\times E\;,
\end{equation}
where boundary values for both fields are given at the tent bottom and
$\nabla=\nabla_x$ denotes the spatial gradient.

The approach of MTP schemes is to map the tent domain to
a spacetime cylinder $\omega_V \times (0,1)$ and solve the
transformed equation there.
The transformation from the cylinder to the tent is denoted by
$\Phi : \omega_V \times (0,1)   \rightarrow  K$ and is defined by
$\Phi(x, \hat t) = (x, \varphi(x,\hat t))$ where
\begin{equation*}
\varphi(x,\hat t) =
(1-\hat t) \varphi_b(x) + \hat t \varphi_t (x)\;.
\end{equation*}
It is similar to the Duffy transformation mapping a square to a
triangle.

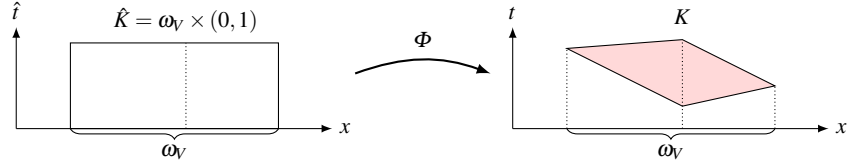
\begin{figure}
\centering
\begin{tikzpicture}[scale=3]
  \tikzstyle{redfill} = [fill=red!30,fill opacity=0.5]

  \begin{scope}[yshift=-3]
  \coordinate (p10) at (5.952,0.205);
  \coordinate (p18) at (5.440,0.461);
  \coordinate (p19) at (6.362,0.295);
  \coordinate (p24) at (5.952,0.500);
  \filldraw[redfill] (p10) -- (p18) -- (p24) -- (p19) -- cycle;
\end{scope}
  \draw[densely dotted] (5.440,0) -- (p18);
  \draw[densely dotted] (5.952,0) -- (p24);
  \draw[densely dotted] (6.362,0) -- (p19);
    
  \draw [decorate,decoration={brace,amplitude=4pt,mirror}] (5.440,0) -- (6.362,0) node[midway,below,yshift=-3pt] {$\omega_V$};
  \node at ($(p24)+(0,0.1)$) {$K$};
  \coordinate (o) at (5.2,0);
  \draw[-latex] (o)--($(o)+(1.4,0)$) node[right] {$x$};
  \draw[-latex] (o)--($(o)+(0,0.45)$) node[above] {$t$};

  \coordinate (b1) at (5.440,0);
  \coordinate (b2) at (5.952,0);
  \coordinate (b3) at (6.362,0);
  \coordinate (t1) at (5.440,0.38);
  \coordinate (t2) at (5.952,0.38);
  \coordinate (t3) at (6.362,0.38);
  \coordinate (shift) at (-2.2,0);
  \coordinate (oref) at ($(o)+(shift)$);
  \draw ($(b1)+(shift)$) -- ($(b2)+(shift)$) -- ($(b3)+(shift)$) --($(t3)+(shift)$) -- ($(t2)+(shift)$) -- ($(t1)+(shift)$) -- cycle;
  \draw[densely dotted] ($(b2)+(shift)$) -- ($(t2)+(shift)$);

  \draw [decorate,decoration={brace,amplitude=4pt,mirror}] ($(b1)+(shift)$) -- ($(b3)+(shift)$) node[midway,below,yshift=-3pt] {$\omega_V$};
  \node at ($(t2)+(shift)+(0,0.1)$) {$\hat K = \omega_V \times (0,1)$};

  \draw[-latex] (oref) -- ($(oref)+(1.4,0)$) node[right] {$x$};
  \draw[-latex] (oref) -- ($(oref)+(0,0.45)$) node[above] {$\hat t$};

  \coordinate (mp) at ($(p24)+0.5*(shift)$);
  \draw[-latex,thick,bend angle=20, bend left] ($(mp)+(-0.35,-0.15)$) to ($(mp)+(0.25,-0.15)$);
  \node (phi) at ($(mp)+(-0.05,0)$) {$\Phi$};
  
\end{tikzpicture}
\caption{Tent mapped from a tensor product domain. \label{mtp_sat:figmap}}
\end{figure}

With the notation
\begin{equation*}
\opskew E = \left( \begin{array}{ccc}
                     0 & E_z & -E_y \\
                     -E_z & 0 & E_x \\
                     E_y & -E_x & 0
\end{array} \right)\;,
\end{equation*}
we can rephrase the curl operator as 
$
\nabla \times E = \opdiv \opskew E,
$
where the divergence of the matrix function is taken row-wise.
To simplify notation further, we define
$u : K \rightarrow
{\setR}^6$ by 
$u = (E, H),$ and set 
$g : K \rightarrow \setR^6$ and $f : K \rightarrow \setR^{6
  \times 3}$ by
\begin{equation}
  \label{mtp_sat:maxwell_def_gf}
  g(u) = 
  \begin{bmatrix}
    \varepsilon E \\ \mu H  
  \end{bmatrix}\;,
  \qquad
  f(u) = \left[ \begin{array}{r} 
                  -\opskew H \\
                  \opskew E 
                \end{array} \right]\;.
\end{equation}
Then~\eqref{mtp_sat:maxwell} may be rewritten as
the conservation law
$\partial_t g(u) + \opdiv_x f(u) = 0.$
Furthermore, we define $F(u) \in \setR^{6 \times4}$ as
\begin{equation*}
F(u) = 
\left[ f(u) \;\; g(u) \right] = 
\left[ \begin{array}{rc} 
    -\opskew H & \;\varepsilon E \\
      \opskew E & \;\mu H 
\end{array} \right]\;,
\end{equation*}
which allows us to write Maxwell's system (\ref{mtp_sat:maxwell}) 
as the spacetime conservation law
\begin{equation}
  \label{mtp_sat:divxt}
\opdiv_{x,t} F(u) = 0\;.
\end{equation}
For each row of $F$, the spacetime divergence $\opdiv_{x,t}$
sums the spatial divergence
of the first three components with the time-derivative of the last
component.

Now, we apply the Piola transformation to pull back $F$ from the tent $K$ to
the cylinder using the mapping $\Phi$. The derivative of $\Phi$ and
its transposed inverse are
\begin{equation*}
\Phi^\prime =  \begin{bmatrix}{} I & 0 \\ 
                        \nabla \varphi^T & \delta \end{bmatrix} 
\qquad \text{and} \qquad
(\Phi^\prime)^{-T} =   \begin{bmatrix} I &  \;\;-\delta^{-1} \, \nabla \varphi \\ 
                               0 & \delta^{-1} \end{bmatrix} \;.
\end{equation*}
The Piola transform of $F$ is
$   \widehat F(\hat u) = {\mathcal P} \{ F \}  = (\det \Phi^\prime) (F
   \circ \Phi) (\Phi^\prime)^{-T}  $
with $\hat u=u\circ \Phi$.
Since the Piola transform provides an algebraic transformation of the
divergence, equation~\eqref{mtp_sat:divxt} is  simply transformed to 
$\opdiv_{x, \hat t} \widehat F(\hat u) = 0$
on the spacetime cylinder.
Then, inserting the Jacobian of $\Phi$ leads us to the transformed equation
\begin{equation} \label{mtp_sat:transformedequ}
\partial_{\hat t} ( g(\hat u) - f(\hat u) \nabla\varphi) +
\opdiv_{x} ( \delta f(\hat u)) = 0\;,
\end{equation}
where  $\delta(x) = \varphi_t(x) - \varphi_b(x)$ is the local height of
the tent. Note that $\nabla\varphi$ is an affine-linear function in
quasi-time $\hat t$.
Equation (\ref{mtp_sat:transformedequ}) describes the evolution of
$\hat u$ along quasi-time from $\hat t = 0$ to $\hat t = 1$.
Details of the calculations are given in \cite{GSW17}.

The next step is the space discretization
of~(\ref{mtp_sat:transformedequ})
by a standard discontinuous Galerkin method.
Let $V_h \subset [L_2]^6$ be the DG finite element space of degree $p$ 
on~$\mathcal{T}$.
On each tent we search for $\hat u : [0,1] \rightarrow V_h$ such that
\begin{equation*}
\int_{\omega_V} \partial_{\hat t} \big[ g(\hat u) - f(\hat u) \nabla \varphi\big] \, v_h -
\sum_{T \subset \omega_{V}} \int_T \delta f(\hat u) \nabla v_h + 
\sum_{F \subset \omega_V} \int_F \delta f_n (\hat u^+, \hat u^-) \llbracket v \rrbracket =0
\end{equation*}
holds for all $v_h \in V_h$ and all $\hat t \in [0,1]$. Only the restriction of 
$V_h$  on the patch $\omega_V$ is used in this equation. The numerical flux
$f_n(\hat u^+, \hat u^-)$ depends on the positive trace $\lim_{s \rightarrow
  0^+} \hat u(x+sn)$ and negative trace $\lim_{s \rightarrow 0^+} \hat u(x-sn)$,
where $n$ is a unit normal vector of arbitrary orientation to the
face.  The jump is defined as usual by
 $\llbracket \hat u \rrbracket := \hat u^+ - \hat u^-$ and the mean value by $\{ \hat u \}
:= \tfrac{1}{2}( \hat u^+ + \hat u^- )$. One example is the upwind 
flux~\cite[p.~434]{HW08}
\begin{equation*}
f_n(\hat u^+, \hat u^-) = \left[ \begin{array} {r}
                         \{ \hat H \} \times n + \llbracket \hat E_t \rrbracket \\
                         -\{ \hat E \} \times n + \llbracket \hat H_t \rrbracket
                         \end{array} \right]\;,
\end{equation*}
with the tangential components $\hat E_t=-(\hat E\times n)\times n$
and $\hat H_t=-(\hat H\times n)\times n$ of $\hat E = E\circ\Phi$ and $\hat H = H\circ\Phi$.
Note that the local tent height $\delta$
enters the boundary integrals as a multiplicative factor. At the outer
boundary of the vertex patch we have $\delta = 0$,
so the facet integrals on the outer boundary disappear.
For the above semidiscrete system,
initial values for the tent problem are given finite element functions at the
tent bottom. The finite element solution on the tent top provides the
initial conditions for the next level tent. Therefore, no projection of initial
values is needed when propagating from one tent to the next.

After the semi-discretization, as usual, we
are left to solve a system of
$N = \dim V_h(\omega_V)$ ordinary differential equations for $U : [0,1]
\rightarrow \setR^N$,
\begin{equation} \label{mtp_sat:equ_ode}
\frac{d}{d\hat t} \left[ M U \right](\hat t) - A U(\hat t) = 0 \;, \qquad
\hat t \in (0,1)\;,
\end{equation}
given $U(0)$. The non-standard feature of~\eqref{mtp_sat:equ_ode}
is that $M$
 is an affine-linear function of 
the  quasi-time $\hat t$ 
(since our mapping enters the mass matrix $M$ through $\nabla \varphi$).
The
matrix $A$ is independent of $\hat t$.
A straightforward approach is to substitute $Y = M U$ and solve
\begin{equation*}
\frac{d}{d \hat t} Y - A M^{-1} Y = 0 \;,
\end{equation*}
instead of~\eqref{mtp_sat:equ_ode}. 
Although first order convergence was observed with
this strategy, further numerical studies showed reduced order of
convergence if the stage-order of the Runge Kutta (RK) method is not high
enough -- see Fig.~\ref{mtp_sat:wave2d_conv} (right).
While the implicit MTP schemes discussed in~\cite{GSW17} do not show
this problem, the issue remains critical for explicit schemes. Thus, we
propose to use a new type of  explicit time-stepping for time discretization,
discussed next.

\section{Structure-aware Taylor time-stepping}

Returning to the ordinary differential equation
(\ref{mtp_sat:equ_ode}) and continuing to make the substitution $Y =
M U$, we now reconsider the previous equation
as the following differential-algebraic system:
\begin{equation}
  \label{mtp_sat:system_ode}
  \frac{d}{d\that} Y = A U\;, \qquad Y =  M U\;.
\end{equation}
We begin by subdividing the interval $(0,1)$ into
$m\in\mathbb{N}$ smaller intervals of size $\frac{1}{m}$, defined
by $(\that_i,\that_{i+1})=(\frac{i}{m},\frac{i+1}{m})$, for
$i\in\mathbb{N}$ and $0\le i\le m-1$.
Recall that  $A$ is independent of quasi-time $\that$,
and $M$ is an affine function  of $\that$,  i.e.,
\[
  M(\hat t) = M_i + (\hat t - \that_i) M', \qquad \hat t \in (\that_i, \that_{i+1})
\]
where $M_i = M(\that_i)$ and the derivative $M'$ is a constant matrix.
We want to design a time-stepping scheme that is aware of this
structure.

Consider the approximations to  $Y, U$ on $(\that_i,\that_{i+1})$
in the form of Taylor polynomials $Y_i, U_i$ of degree $q$, defined by
\begin{equation}
  \label{mtp_sat:taylor_pol}
  Y_i(\hat t) = \sum_{n=0}^{q} \frac{(\hat t - \hat t_i)^n}{n!} Y_{i,n} 
  \qquad 
  U_i(\that) =  \sum_{n=0}^{q-1} \frac{(\hat t - \hat t_i)^n}{n!}
  U_{i,n}\;,
  \qquad \hat t \in (\that_i,\that_{i+1})\;,
\end{equation}
where
$Y_{i,n}=Y_i^{(n)}(\that_i)$ and $U_{i,n}=U_i^{(n)}(\that_i).$ To find
these derivatives, we differentiate both equations
of (\ref{mtp_sat:system_ode}) $n$ times to get 
\begin{align*}
Y^{(n+1)}(\that) & = A U^{(n)}(\that)\;, &n \ge 0\;, \\
Y^{(n)}(\that) & = M(\that) U^{(n)}(\that) + n M' U^{(n-1)}(\that)\;, &n \ge 1\;.
\end{align*}
For the second equation we used Leibnitz' formula
$(fg)^{(n)} = \sum_{i=0}^n { n \choose i } f^{(i)} g^{(n-i)}$, and the
fact that $M$ is affine-linear.  Evaluating these  equations
for the Taylor polynomials $Y_i, U_i$ at $\that = \that_i$,
we obtain a recursive formula for $Y_{i, n}$ and $U_{i, n}$ in terms
of $U_{i, n-1}$, namely
\begin{equation}
  \label{mtp_sat:coeff_rec}
    \begin{array}{r c l @{\qquad}l}
      Y_{i,n} & = & A U_{i,n-1}\;, & 1\le n \le q\;, \\
      M_i U_{i,n} & = & Y_{i,n} - n M' U_{i,n-1}\;, & 1\le n \le q-1\;,
    \end{array}
\end{equation}
for all $0\le i\le m-1$. Given $Y_{0,0} = Y(\hat t_0)$, $M_0U_{0,0} =
Y_{0,0}$, applying~\eqref{mtp_sat:coeff_rec} with $i=0$ gives the
approximate functions $Y_0(\hat t), U_0(\hat t)$ in the first subinterval
$(\that_0, \that_1)$. 
The recursive formulas are  initiated for later subintervals 
at $n=0$ by 
\begin{align} \label{mtp_sat:recursion-init}
  Y_{i,0} &= Y_{i-1}(\that_i), 
  &
    M_i U_{i,0} &  =  Y_{i,0}\;, && 1\le i\le m-1\;.
\end{align}
After the final subinterval, we get
$Y_{m-1}(t_m)$, our approximation to $Y(1)$.
We shall refer to 
the new time-stepping  scheme generated
by~\eqref{mtp_sat:coeff_rec}
as the $q$-stage {\em SAT (structure-aware Taylor) time-stepping.} 

Note that $Y_{m-1}(t_m)$ is our approximation to $Y = MU$ at the 
top of the
tent. This value is then passed to the next tent in time. The time
dependence of $M$ arises from the time dependence of
$\nabla\varphi$. This gradient is continuous along spacetime lines of constant
spatial coordinates. Therefore, when passing from one
element of a tent to the same element within the next tent in
time,  $Y$ is continuous (since the solution $U$ is continuous).
Of course, on
flat fronts $\nabla \varphi = \nabla \tau = 0$, so there $M$ is just a
diagonal matrix containing the material parameters.

To briefly remark on the expected convergence rate of a
$q$-stage SAT time-stepping, recall that  due to the mapping of the MTP method we
solve for $\hat u = u \circ \Phi$, which satisfies
$  \partial^n_{\that} \hat u = \delta^n (\partial^n_tu) \circ \Phi$.
The causality condition implies that $\delta \rightarrow 0$ if the
mesh size $h \rightarrow 0$. Thus we may expect the 
$n^{th}$ temporal derivative of $\hat u$, and correspondingly $U^{(n)},$ to 
go to zero at the rate $\mathcal{O}(h^n)$.
By using a $q$-stage SAT time-stepping, we approximate the
first $q-1$ terms of the exact Taylor expansion of $U$. Thus we expect the
convergence rate to be $O(h^q)$, the size of the remainder term involving
$U^{(q)}$. The next  section
provides numerical evidence for this.

Before concluding this section, we should note that in~\eqref{mtp_sat:coeff_rec} and~\eqref{mtp_sat:recursion-init}, we
tacitly assumed that $M_i$ is invertible. Let us show that this is
indeed the case whenever the causality condition (see
\S\ref{sec:mesh}) $|\nabla\varphi| < \sqrt{\varepsilon\mu}$ is
fulfilled. At any quasi-time $\hat t$,
given a $\hat w = (\hat{w}_E, \hat{w}_H) \in V_h$ whose coefficient vector in the basis
expansion is $W \in \mathbb{R}^N$,  consider the equation  
$M(\hat t) U = W$ for the coefficient vector $U$ of $ \hat u \in V_h$.
This equation, in variational form, is
\begin{equation}
  \label{mtp_sat:a}
\int_{\omega_V} [ g(\uhat) - f(\uhat) \nabla
\varphi ] \cdot \hat v=
\int_{\omega_V} (\hat{w}_E, \hat{w}_H) \cdot \hat v,
\qquad \text{ for all } \hat v \in V_h.
\end{equation}
Let 
$a(\hat u, \hat v)$ denote the left hand side of~\eqref{mtp_sat:a}.
To prove solvability of~\eqref{mtp_sat:a}, it suffices
to prove that $a(\cdot,\cdot)$ is a coercive bilinear form on
 $[L_2]^6$ for any $\hat t$.
By inserting $g(\uhat) = [\varepsilon \Ehat, \mu \Hhat ]^T$ and
$f(\uhat) = [-\opskew \Hhat, \opskew \Ehat]^T$ into $a(\uhat, \uhat)$,
\begin{eqnarray*}
  a(\uhat,\uhat) & = & \int_{\omega_V}  ( \varepsilon \Ehat - \Hhat \times \nabla \varphi)\cdot \Ehat +
                       ( \mu \Hhat + \Ehat \times \nabla \varphi)\cdot \Hhat \\
                 & = & \int_{\omega_V}  \varepsilon \Ehat\cdot \Ehat + \mu \Hhat\cdot \Hhat + 2 (\Ehat \times \nabla \varphi)\cdot \Hhat \\
                 & \ge & \int_{\omega_V}  \varepsilon \Ehat\cdot \Ehat + \mu \Hhat\cdot \Hhat
                         - 2 \frac{|\nabla \varphi|}{\sqrt{\varepsilon\mu}} \sqrt{\varepsilon}|\Ehat|\sqrt{\mu}|\Hhat| \;,
\end{eqnarray*}
where we used the Cauchy-Schwarz inequality and inserted
$\sqrt{\varepsilon}$ and $\sqrt{\mu}$ to achieve the desired
scaling. By applying Young's inequality and 
$|\nabla\varphi| < \sqrt{\varepsilon\mu}$, 
\begin{eqnarray*}
  a(\uhat,\uhat) & \ge & \int_{\omega_V}  \varepsilon \Ehat\cdot \Ehat + \mu \Hhat\cdot \Hhat
                         - \frac{|\nabla \varphi|}{\sqrt{\varepsilon\mu}} (\varepsilon \Ehat\cdot \Ehat + \mu \Hhat\cdot \Hhat) \\
                 & = & \int_{\omega_V}  \left( 1 - \frac{|\nabla \varphi|}{\sqrt{\varepsilon\mu}}\right) (\varepsilon \Ehat\cdot \Ehat + \mu \Hhat\cdot \Hhat)
                       \ge C \min{(\varepsilon,\mu)} \|\uhat\|_{L_2}^2\;,
\end{eqnarray*}
form some constant  $C>0$.
Thus $M_i$ is invertible and the SAT time-stepping 
is well defined on all tents respecting the causality
condition.

One may exploit the specific
details of the Maxwell problem 
to avoid the assembly and the inversion of matrices $M_i$ (as we have
done in our implementation). In fact,
instead of~\eqref{mtp_sat:a}, we can explicitly solve the corresponding exact
undiscretized equation obtained by replacing $V_h$ by $[L_2]^6$
in~\eqref{mtp_sat:a}. The solution $\hat u = ( \hat{E}, \hat{H})$ in closed
form reads
\begin{eqnarray*}
  \Ehat & = & \frac{1}{\varepsilon \mu - | \nabla \varphi |^2 } \left(  I -
              \frac{1}{\varepsilon \mu} \nabla \varphi \nabla \varphi^T
              \right)   (\mu \hat w_E + \hat w_H \times \nabla \varphi) \;, \\
  \Hhat & = & \frac{1}{\varepsilon \mu - | \nabla \varphi |^2 } \left(  I -
              \frac{1}{\varepsilon \mu} \nabla \varphi \nabla \varphi^T
              \right)   (\varepsilon \hat w_H - \hat w_E \times \nabla \varphi)\;.
\end{eqnarray*}
We then perform a projection of these into $V_h$
to obtain the coefficients $U(\that_i)$.
For uncurved elements, this just involves the inversion of a diagonal
mass matrix. For the small number of curved elements, we use a
highly optimized algorithm which uses an approximation instead of the
exact inverse mass matrix.

\section{Numerical Results}
The MTP discretization in combination with the SAT 
time-stepping on tents is implemented within the Netgen/NGSolve finite
element library. In this section numerical results concerning accuracy
as well as performance are reported.

\subsection{Convergence studies in two space dimensions}
We consider the model problem in two space dimensions
\begin{equation*}
\partial_t \varepsilon E_z = \partial_xH_{y} - \partial_yH_{x}\;, \qquad 
\partial_t \mu H_x = -\partial_y E_{z}\;,  \qquad
\partial_t \mu H_y = \partial_x E_{z}\;,
\end{equation*}
on the spacetime cube $[0,\pi]^2 \times [0,\sqrt{2}\pi]$.
Parameters are set $\varepsilon = \mu = 1$ such
that speed of light is $c = 1$.
Initial and boundary values are set such that the exact solution
is given by
\begin{eqnarray*}
E_z & =  & \sin(x) \sin(y) \cos{(\sqrt{2}t)}\;, \\
H_x & =  & -\tfrac{1}{\sqrt{2}} \sin(x) \cos(y) \sin{(\sqrt{2}t)}\;, \\
H_y & =  & \tfrac{1}{\sqrt{2}} \cos(x) \sin(y) \sin{(\sqrt{2}t)}\;.
\end{eqnarray*}

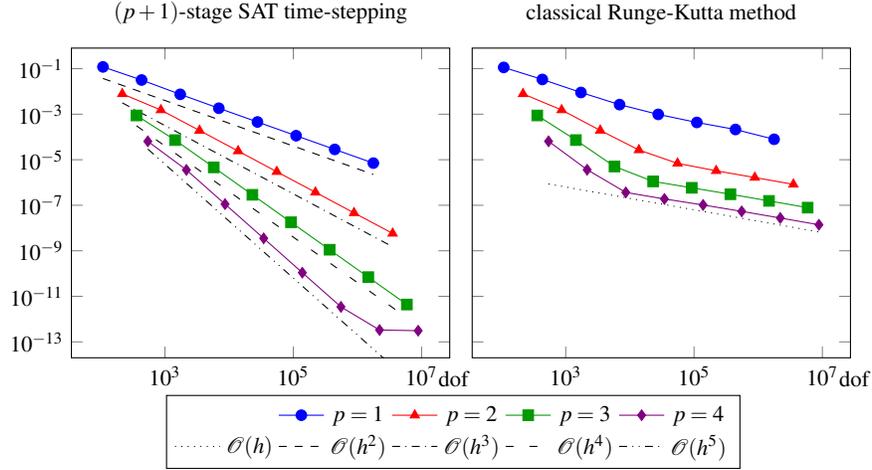
\begin{figure}
  \centering
  \begin{tikzpicture}
    \begin{groupplot}[%
      group style={%
        group name={},%
        group size=2 by 1,%
        horizontal sep=1em,%
        yticklabels at=edge left,%
      },%
      width=0.565\textwidth,%
      ]%
      \nextgroupplot[xmode=log,ymode=log,%
      tick label style={font=\footnotesize},%
      xtick={1e3,1e5,1e7},%
      ymax=0.8, ymin=2e-14,%
      ytick={1e-1,1e-3,1e-5,1e-7,1e-9,1e-11,1e-13},%
      legend style = {at={(0.25,-0.12)},anchor=north west},%
      legend columns = 5,
      title={$(p+1)$-stage SAT time-stepping},
      xlabel={\footnotesize dof},
      xlabel style = {at={(1.01,0.01)} ,anchor=south},
      ]%
      \addlegendimage{empty legend}
      \addlegendentry{}
      \addplot[color=blue,mark=*] table[x=p1_ndof,y=p1_l2error] {table_err_wave2d_taylor_1.dat};
      \addlegendentry{$p=1$};
      \addplot[color=red,mark=triangle*] table[x=p2_ndof,y=p2_l2error] {table_err_wave2d_taylor_1.dat};
      \addlegendentry{$p=2$};
      \addplot[color=green!60!black,mark=square*] table[x=p3_ndof,y=p3_l2error] {table_err_wave2d_taylor_1.dat};
      \addlegendentry{$p=3$};
      \addplot[color=violet,mark=diamond*] table[x=p4_ndof,y=p4_l2error] {table_err_wave2d_taylor_1.dat};
      \addlegendentry{$p=4$};
      
      \addplot[draw=none,dotted] table[x=p1_ndof,y expr={4/\thisrow{p1_ndof}}] {table_err_wave2d_taylor_1.dat};
      \addlegendentry{$\mathcal{O}(h)$};
      \addplot[dashed] table[x=p1_ndof,y expr={4/\thisrow{p1_ndof}}] {table_err_wave2d_taylor_1.dat};
      \addlegendentry{$\mathcal{O}(h^2)$};
      \addplot[dashdotted] table[x=p2_ndof,y expr={10/\thisrow{p2_ndof}^(3/2)}] {table_err_wave2d_taylor_1.dat};
      \addlegendentry{$\mathcal{O}(h^3)$};
      \addplot[loosely dashed] table[x=p3_ndof,y expr={40/\thisrow{p3_ndof}^2}] {table_err_wave2d_taylor_1.dat};
      \addlegendentry{$\mathcal{O}(h^4)$};
      \addplot[dashdotdotted] table[x=p4_ndof,y expr={2e2/\thisrow{p4_ndof}^(5/2)}] {table_err_wave2d_taylor_1.dat};
      \addlegendentry{$\mathcal{O}(h^5)$};
      \nextgroupplot[xmode=log,ymode=log,%
      tick label style={font=\footnotesize},%
      xtick={1e3,1e5,1e7},%
      ymax=0.8, ymin=2e-14,%
      ytick={1e-1,1e-3,1e-5,1e-7,1e-9,1e-11,1e-13},%
      title=classical Runge-Kutta method,
      xlabel={\footnotesize dof},
      xlabel style = {at={(1.01,0.01)} ,anchor=south},
      ]%
      \addplot[color=blue,mark=*]table[x=p1_ndof,y=p1_l2error] {table_err_wave2d_rk4.dat};
      \addplot[color=red,mark=triangle*]table[x=p2_ndof,y=p2_l2error] {table_err_wave2d_rk4.dat};
      \addplot[color=green!60!black,mark=square*] table[x=p3_ndof,y=p3_l2error] {table_err_wave2d_rk4.dat};
      \addplot[color=violet,mark=diamond*] table[x=p4_ndof,y=p4_l2error] {table_err_wave2d_rk4.dat};
      \addplot[dotted] table[x=p4_ndof,y expr={2e-5/\thisrow{p4_ndof}^(1/2)}] {table_err_wave2d_rk4.dat};
    \end{groupplot}
  \end{tikzpicture}
\caption{Spatial $L_2$ error of all field components over degrees of
  freedom (dof) for the $(p+1)$-stage SAT time-stepping (left) and the classical Runge-Kutta (right).} \label{mtp_sat:wave2d_conv}
\end{figure}

Based on a spatial mesh with mesh size $h$,
we generate a tent pitched mesh such that the maximal slope
$|\nabla \varphi|$ is bounded by $(2c)^{-1}$ and apply a discontinuous
Galerkin method in space using polynomials of order $p$, with
$1 \leq p \leq 4$. On each cylinder we perform a $(p+1)$-stage
SAT time-stepping with $m = 2p$ intervals. The spatial $L_2$ error of all
field components at the final time is reported in the left plot of
Fig.~\ref{mtp_sat:wave2d_conv}. We observe that the error goes to
zero at the optimal rate of $\mathcal{O}(h^{p+1})$ until we are
close to 
machine precision.

In contrast, the right plot in Fig.~\ref{mtp_sat:wave2d_conv}
illustrates the previously mentioned loss of convergence rates
when the classical Runge-Kutta method is used. The convergence rates
stagnate at first order no matter what $p$ is used. A similar  behavior was
also observed for other explicit Runge-Kutta methods.

\subsection{Large scale problem in three space dimensions}
\begin{figure}
  \centering
  \begin{tikzpicture}[scale=1]
    \node[circle,draw, inner sep=28, semithick,
    path picture={
      \node at ($(path picture bounding box.center)+(-0.6,0.0)$) {
        \includegraphics[width=.55\textwidth]{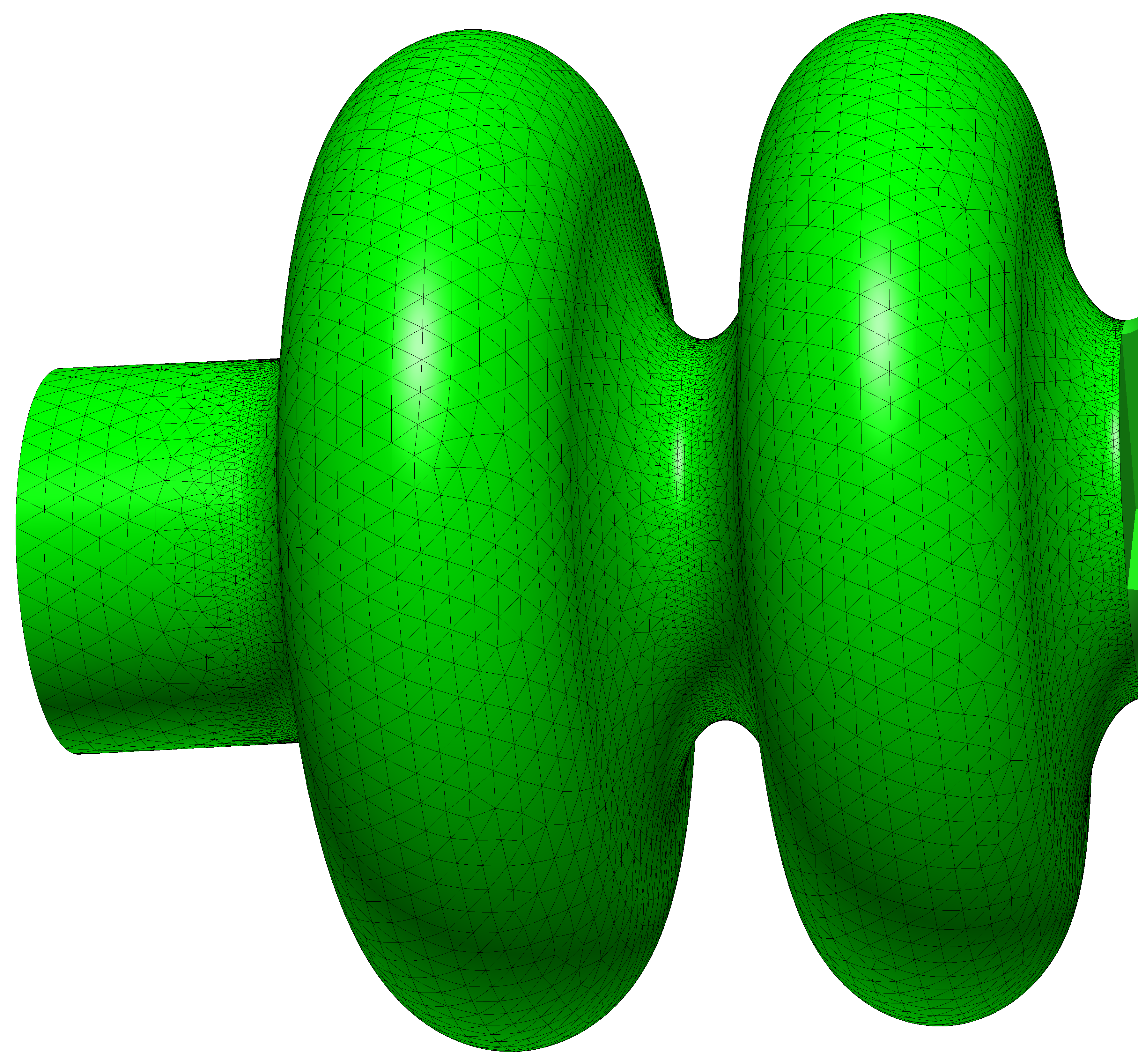}
      };
    }] (A) {};
    \node at (5.8,0) {\includegraphics[width=0.75\textwidth, trim=0 0 100 0, clip]{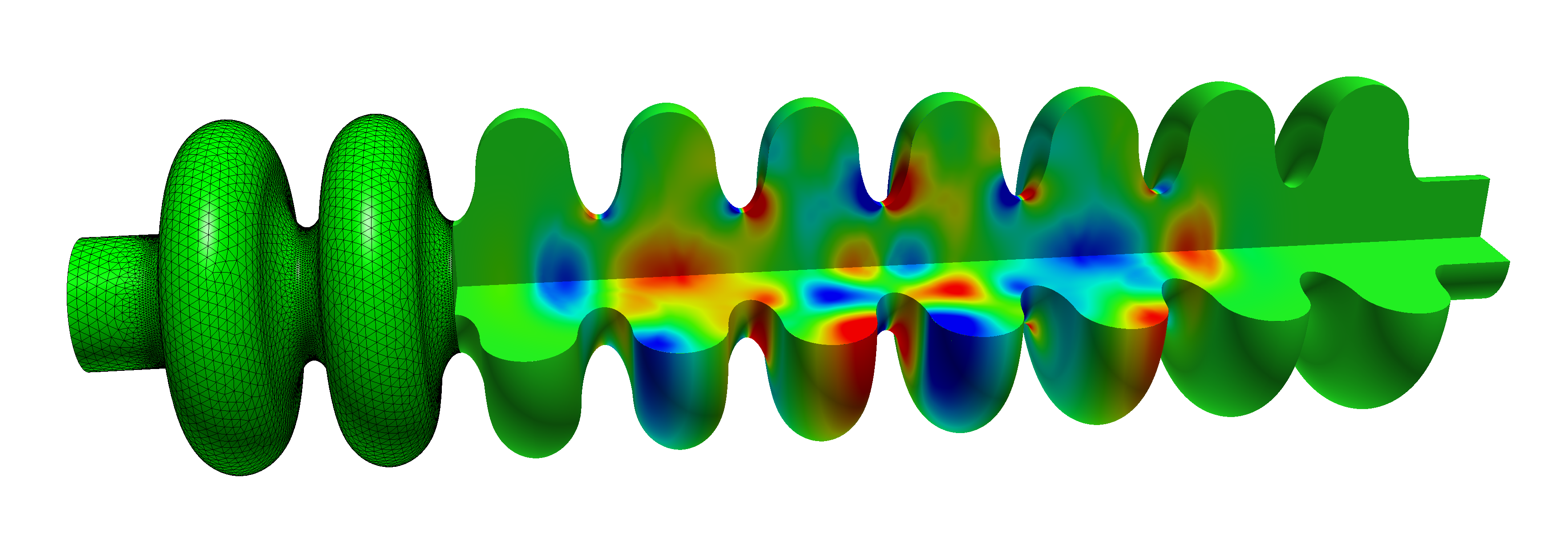}};
  \end{tikzpicture}
  \caption{Tetrahedral mesh with 489k curved elements, ratio of the largest to the smallest element of approximately
    5:1 and the $H_y$ component of solution at $t=260$ calculated with spatial polynomial order $p=3$.}
  \label{mtp_sat:resonator_smooth}
\end{figure}

As a second example we present a simulation on a domain similar to the
resonator shown in \cite{HPSTW15}. The geometry is given as body of
revolution of smooth B-spline curves. The mesh consisting of 489593
curved tetrahedral elements is shown in
Fig.~\ref{mtp_sat:resonator_smooth}. Due to higher curvature the mesh is
refined along the inner roundings, where the ratio of the largest to
the smallest element is approximately 5:1.  We used
a Gaussian peak
(located at the axis of
revolution and the position of the fifth inner rounding)
for the electric field as initial data. The explicit MTP 
scheme with SAT time-stepping then computed the solution
at $t=260$ using time slabs of height 1,
with each slab composed of $N_{\mathrm{tents}}=149072$ tents. On each
tent we used a $(p+1)$-stage SAT time-stepping with $m=2p$ intervals,
where $p$ denotes the spatial polynomial order.
With the spatial degrees of freedom $N_{\mathrm{dof},i}$ of
the $i^{th}$ tent and the number of stages $q=p+1$,
we obtain the total spacetime degrees of freedom per time slab
\begin{equation*}
  \sum_{i=1}^{N_{\mathrm{tents}}} N_{\mathrm{dof},i}\,m\, q =
  \left(\sum_{i=1}^{N_{\mathrm{tents}}} N_{\mathrm{dof},i}\right) 2p(p+1)\;.
\end{equation*}
The corresponding numbers of degrees of freedom and the simulation times
are shown in Table \ref{mtp_sat:table_smooth}. In \cite{HPSTW15}
a similar problem is solved  using a discontinuous Galerkin method
with quadratic elements, combined with a polynomial Krylov subspace
method in time. Using 96 cores it took them 7:10 hours to reach the final
time. Our simulation with polynomial order $p=3$, which has a
comparable number of unknowns, took 3:33 hours on 64 cores. This 
significant speed up is an illustration of  the capability of the new method.
The $H_y$ component of the obtained solution at $t=260$, using third order
polynomials in space, is shown in Fig.~\ref{mtp_sat:resonator_smooth}.
\begin{table}[h]
  \centering
  \begin{tabular}{@{~~}l@{\quad} |@{\quad}l@{\quad} |@{\quad}l@{~~}}
    & $p = 2$ & $p = 3$ \\
    \hline\hline
    number of spatial dof & $2.938\times 10^7$ &  $5.875\times 10^7$ \\
    number of spacetime dof per slab & $1.908\times 10^9$ & $7.632\times 10^9$ \\
    simulation time per slab & 4.6 s & 49.2 s\\
    total simulation time & 20 min & 3 h 33 min
  \end{tabular}
  \caption{Number of degrees of freedom and simulation times for
    spatial polynomial orders $p=2,3$. This data was generated using a
    shared memory server with 4 E7-8867 CPUs with 16 cores each.}
  \label{mtp_sat:table_smooth}
\end{table}
\vspace{-30pt}

\vspace{-0.5cm}
\begin{acknowledgement}
  This work was supported in part by the National Science Foundation.
\end{acknowledgement}
\vspace{-0.6cm}



\end{document}